\documentclass[11pt,reqno]{amsart}
\usepackage{amsfonts,amssymb,amsmath}

\newtheorem{theorem}{Theorem}
\newtheorem{corollary}[theorem]{Corollary}
\newtheorem{lemma}[theorem]{Lemma}
\newtheorem{proposition}[theorem]{Proposition}
\newtheorem{remark}[theorem]{Remark}

\numberwithin{theorem}{section}

\begin{document}

\title{Neighboring Fractions in Farey Subsequences}

\author{Andrey O. Matveev}
\address{Data-Center Co., RU-620034, P.O.~Box~5, Ekaterinburg,
Russian~Federation} \email{aomatveev@dc.ru aomatveev@hotmail.com}


\begin{abstract}
We present explicit formulas for the computation of the~neighbors of several elements of Farey subsequences.
\end{abstract}

\maketitle

\section{Introduction}

The~{\em Farey sequence $\mathcal{F}_m$ of order $m$\/} is the~ascending set of rational numbers~$\tfrac{h}{k}$, written in reduced terms, such that $\tfrac{0}{1}\leq\tfrac{h}{k}\leq\tfrac{1}{1}$ and $1\leq
k\leq m$, see,~e.g.,~\cite[Chapter~27]{B}, \cite{CZ}, \cite[\S{}3]{GCh}, \cite[Chapter~4]{GKP}, \cite[Chapter~III]{HW}, \cite[Chapter~6]{NZM}, \cite[Chapter~6]{P}, \cite[Sequences A006842 and A006843]{S}, \cite[Chapter~5]{T}. For example,
\begin{equation*}
\mathcal{F}_5=\left(\tfrac{0}{1}<\tfrac{1}{5}<\tfrac{1}{4}<\tfrac{1}{3}<\tfrac{2}{5}<\tfrac{1}{2}
<\tfrac{3}{5}<\tfrac{2}{3}<\tfrac{3}{4}<\tfrac{4}{5}<\tfrac{1}{1}\right)\ .
\end{equation*}

Recall that the~map
\begin{equation}
\label{eq:10} \mathcal{F}_m\to\mathcal{F}_m\ ,\ \ \  \tfrac{h}{k}\mapsto\tfrac{k-h}{k}\ ,
\end{equation}
is order-reversing and bijective.

Let $\overline{\mu}(\cdot)$ denote the~M\"{o}bius function on positive integers. We use the~notation $[s,t]$
to denote the~interval $\{s,s+1,\ldots,t\}$ of positive integers; the~greatest common divisor of $s$ and $t$ is denoted by $\gcd(s,t)$, and we write $s|t$ if $t$ is divisible by $s$.

Notice that if $\mathrm{y}$ is a formal variable, then for a positive integer $i$ and for a~nonempty interval $[i'+1,i'']$ we have
\begin{equation*}
\sum_{\substack{j\in[i'+1,i'']:\\ \gcd(i,j)=1}}\mathrm{y}^j=
\sum_{\substack{d\in[1,i]:\\
d|i}}\overline{\mu}(d)\frac{\mathrm{y}^{d\left\lceil\frac{i'+1}{d}\right\rceil}-\mathrm{y}^{d\cdot(\left\lfloor
\frac{i''}{d}\right\rfloor+1)}}{1-\mathrm{y}^d}\ .
\end{equation*}
Recall that $\bigl|\{j\in[i'+1,i'']:\ \gcd(i,j)=1\}\bigr|=\sum_{d\in[1,i]:\  d|i}\overline{\mu}(d)
\Bigl(\left\lfloor\frac{i''}{d}\right\rfloor$ $-\left\lfloor\frac{i'}{d}\right\rfloor\Bigr)$,
see, e.g.,~\cite{AM-R}, \cite{AM-I}; see also, e.g.,~\cite{A-K}, \cite{ElB-S} and references therein
on relatively prime sets of integers and, in particular,
on enumeration of $l$-subsets of $[i'+1,i'']$ that are relatively prime to $i$, for any $l$.

If~$\mathrm{x}$ and $\mathrm{y}$ are
formal variables, then we have
\begin{align*}
\sum_{\substack{\frac{h}{k}\in\mathcal{F}_m:\\ \frac{0}{1}<\frac{h}{k}<\frac{1}{2}}}\mathrm{x}^h\mathrm{y}^k
&=\sum_{i\in[1,\lceil\frac{m}{2}\rceil-1]}\mathrm{x}^i\sum_{\substack{j\in[2i+1,m]:\\
\gcd(i,j)=1}}\mathrm{y}^j\\&=\sum_{i\in[1,\lceil\frac{m}{2}\rceil-1]}\mathrm{x}^i
\sum_{\substack{d\in[1,i]:\\
d|i}}\overline{\mu}(d)\frac{\mathrm{y}^{d\left\lceil\frac{2i+1}{d}\right\rceil}-\mathrm{y}^{d\cdot(\left\lfloor
\frac{m}{d}\right\rfloor+1)}}{1-\mathrm{y}^d}\ ;
\end{align*}
thus,
\begin{multline*}
\sum_{\substack{\frac{h}{k}\in\mathcal{F}_m:\\ \frac{0}{1}<\frac{h}{k}<\frac{1}{2}}}\mathrm{x}^h\mathrm{y}^k
=\sum_{d\in[1,\lceil\frac{m}{2}\rceil-1]}\overline{\mu}(d)\frac{\mathrm{y}^d}{1-\mathrm{y}^d}
\Biggl(\frac{\mathrm{x}^d\mathrm{y}^{2d}-\mathrm{x}^{d\left\lceil\frac{m}{2d}\right\rceil}
\mathrm{y}^{2d\left\lceil\frac{m}{2d}\right\rceil}}{1-\mathrm{x}^d\mathrm{y}^{2d}}\\ -\frac{\mathrm{x}^d
-\mathrm{x}^{d\left\lceil\frac{m}{2d}\right\rceil}}{1-\mathrm{x}^d}
\mathrm{y}^{d\left\lfloor\frac{m}{d}\right\rfloor}\Biggr)\ ;
\end{multline*}
for example,
\begin{equation*}
\sum_{\substack{\frac{h}{k}\in\mathcal{F}_5:\\ \frac{0}{1}<\frac{h}{k}<\frac{1}{2}}}\mathrm{x}^h\mathrm{y}^k=\mathrm{x}\mathrm{y}^5+
\mathrm{x}\mathrm{y}^4+\mathrm{x}\mathrm{y}^3+\mathrm{x}^2\mathrm{y}^5=
\frac{\mathrm{y}}{1-\mathrm{y}}
\Biggl(\frac{\mathrm{x}\mathrm{y}^2-\mathrm{x}^3\mathrm{y}^6}{1-\mathrm{x}\mathrm{y}^2}-\frac{\mathrm{x}
-\mathrm{x}^3}{1-\mathrm{x}}\mathrm{y}^5\Biggr)\ .
\end{equation*}
Similarly,
\begin{multline*}
\sum_{\substack{\frac{h}{k}\in\mathcal{F}_m:\\ \frac{1}{2}<\frac{h}{k}<\frac{1}{1}}}\mathrm{x}^h\mathrm{y}^k
=\sum_{d\in[1,\lceil\frac{m}{2}\rceil-1]}\overline{\mu}(d)\frac{\mathrm{x}^d\mathrm{y}^d}{1
-\mathrm{x}^d\mathrm{y}^d}
\Biggl(\frac{\mathrm{x}^d\mathrm{y}^{2d}-\mathrm{x}^{d\left\lceil\frac{m}{2d}\right\rceil}
\mathrm{y}^{2d\left\lceil\frac{m}{2d}\right\rceil}}{1-\mathrm{x}^d\mathrm{y}^{2d}}\\ -
\frac{\mathrm{x}^{d\cdot (\left\lfloor\frac{m}{d}\right\rfloor-\left\lceil\frac{m}{2d}\right\rceil+1)}-
\mathrm{x}^{d\left\lfloor\frac{m}{d}\right\rfloor}}{1-\mathrm{x}^d}
\mathrm{y}^{d\left\lfloor\frac{m}{d}\right\rfloor}\Biggr)\ .
\end{multline*}

Standard tools for construction of the~Farey sequences are numerical and matrix recurrences, as well as
computational tree-like structures, see~\cite[Chapter~4]{GKP}.

The~Farey sequence of order $2m$  contains the~subsequence (see~\cite{AM-R,AM-I,AM-II,AM-C-II})
\begin{equation}
\label{eq:22} \mathcal{F}\bigl(\mathbb{B}(2m),m\bigr):=\left(\tfrac{h}{k}\in\mathcal{F}_{2m}:\ k-m\leq h\leq
m \right)\ ,
\end{equation}
which in a~sense provides us with a~dual description of the~Farey sequence of order $m$, and vice versa; the~notation $\mathbb{B}(2m)$ in definition~(\ref{eq:22}) makes reference to the~Boolean lattice of rank $2m$.
Define the~{\em left\/} and {\em right halfsequences\/}~$\mathcal{F}^{\leq\frac{1}{2}}\bigl(\mathbb{B}(2m),m\bigr)$ and
$\mathcal{F}^{\geq\frac{1}{2}}\bigl(\mathbb{B}(2m),m\bigr)$ of sequence~(\ref{eq:22}) by
\begin{align*}
\mathcal{F}^{\leq\frac{1}{2}}\bigl(\mathbb{B}(2m),m\bigr):&=
\left(\tfrac{h}{k}\in\mathcal{F}\bigl(\mathbb{B}(2m),m\bigr):\
\tfrac{h}{k}\leq\tfrac{1}{2}\right)\intertext{and}
\mathcal{F}^{\geq\frac{1}{2}}\bigl(\mathbb{B}(2m),m\bigr):&=\left(\tfrac{h}{k}
\in\mathcal{F}\bigl(\mathbb{B}(2m),m\bigr):\ \tfrac{h}{k}\geq\tfrac{1}{2}\right)\ ,
\end{align*}
respectively. The~numerators of the~fractions of $\mathcal{F}_m$ are the~numerators of the~fractions of
$\mathcal{F}^{\leq\frac{1}{2}}\bigl(\mathbb{B}(2m),m\bigr)$, while the~denominators of the~fractions of~$\mathcal{F}_m$ are the~numerators of the~fractions of
$\mathcal{F}^{\geq\frac{1}{2}}\bigl(\mathbb{B}(2m),m\bigr)$. More precisely, the~following can be
said~\cite{AM-I,AM-C-II}:

\begin{lemma}
\label{prop:5}
The~maps
\begin{align}
\label{eq:11} \mathcal{F}^{\leq\frac{1}{2}}\!\bigl(\mathbb{B}(2m),m\bigr)&\to\mathcal{F}_m\ , &
\tfrac{h}{k}&\mapsto\tfrac{h}{k-h}\ ,\\ \label{eq:12}
\mathcal{F}_m&\to\mathcal{F}^{\leq\frac{1}{2}}\!\bigl(\mathbb{B}(2m),m\bigr)\ , &
\tfrac{h}{k}&\mapsto\tfrac{h}{k+h}\ ,\\ \label{eq:16}
\mathcal{F}^{\geq\frac{1}{2}}\!\bigl(\mathbb{B}(2m),m\bigr)&\to\mathcal{F}_m\ , &
\tfrac{h}{k}&\mapsto\tfrac{2h-k}{h}\ ,\intertext{and} \label{eq:20}
\mathcal{F}_m&\to\mathcal{F}^{\geq\frac{1}{2}}\!\bigl(\mathbb{B}(2m),m\bigr)\ , &
\tfrac{h}{k}&\mapsto\tfrac{k}{2k-h}\ ,
\end{align}
are order-preserving and bijective.

The~maps
\begin{align} \label{eq:15}
\mathcal{F}^{\leq\frac{1}{2}}\!\bigl(\mathbb{B}(2m),m\bigr)&\to\mathcal{F}_m\ , &
\tfrac{h}{k}&\mapsto\tfrac{k-2h}{k-h}\ ,\\ \label{eq:21}
\mathcal{F}_m&\to\mathcal{F}^{\leq\frac{1}{2}}\!\bigl(\mathbb{B}(2m),m\bigr)\ , &
\tfrac{h}{k}&\mapsto\tfrac{k-h}{2k-h}\ ,\\ \label{eq:13}
\mathcal{F}^{\geq\frac{1}{2}}\!\bigl(\mathbb{B}(2m),m\bigr)&\to\mathcal{F}_m\ , &
\tfrac{h}{k}&\mapsto\tfrac{k-h}{h}\ ,\intertext{and} \label{eq:14}
\mathcal{F}_m&\to\mathcal{F}^{\geq\frac{1}{2}}\!\bigl(\mathbb{B}(2m),m\bigr)\ , &
\tfrac{h}{k}&\mapsto\tfrac{k}{k+h}\ ,
\end{align}
are order-reversing and bijective.
\end{lemma}

For example, the~sequence
\begin{multline*}
\mathcal{F}\bigl(\mathbb{B}(10),5\bigr)=\bigl(\tfrac{0}{1}<\tfrac{1}{6}<\tfrac{1}{5}<\tfrac{1}{4}<\tfrac{2}{7}
<\tfrac{1}{3}<\tfrac{3}{8}<\tfrac{2}{5}<\tfrac{3}{7}<\tfrac{4}{9}\\<\tfrac{1}{2}<\tfrac{5}{9}<\tfrac{4}{7}
<\tfrac{3}{5}<\tfrac{5}{8}<\tfrac{2}{3}<\tfrac{5}{7}<\tfrac{3}{4}<\tfrac{4}{5}<\tfrac{5}{6}<\tfrac{1}{1}\bigr)
\end{multline*}
provides us with a~dual description of the~sequence $\mathcal{F}_5$, and vice versa.

Notice that the~map
\begin{equation*}
\mathcal{F}\bigl(\mathbb{B}(2m),m\bigr)\to\mathcal{F}\bigl(\mathbb{B}(2m),m\bigr)\ ,\ \ \
\tfrac{h}{k}\mapsto\tfrac{k-h}{k}\ ,
\end{equation*}
is order-reversing and bijective.

Applications require a~variety of pairs of adjacent fractions within the~Farey (sub)sequences for the~recurrent computation of other fractions to be performed efficiently; we are particularly interested in
a~family of pairs of adjacent fractions within $\mathcal{F}\bigl(\mathbb{B}(2m),m\bigr)$ which can be described more or less explicitly. We make use of machinery of elementary number theory; the~approach is unified, and it consists in the~transfer of the~results of calculations concerning $\mathcal{F}_m$ to
$\mathcal{F}\bigl(\mathbb{B}(2m),m\bigr)$ by means of some monotone bijections collected
in~Lemma~\ref{prop:5}.

In Section~\ref{section:applications} of the~paper, we show how Farey (sub)sequences arise in analysis of collective decision-making procedures.

In Section~\ref{section:1}, we first recall formulas for the~computation of the~neighbors of
elements of $\mathcal{F}_m$ (Lemma~\ref{prop:1}); then, we present formulas that describe the~neighbors of
fractions of the~form $\tfrac{1}{j}$, $\tfrac{j-1}{j}$ (Corollary~\ref{prop:2}), and of the~form~$\tfrac{2}{j}$, $\tfrac{j-2}{j}$ (Corollary~\ref{prop:7}).

In\hfill Section~\ref{section:2},\hfill we\hfill present\hfill formulas\hfill describing\hfill consecutive fractions
\\ in~$\mathcal{F}\bigl(\mathbb{B}(2m),m\bigr)$ (Proposition~\ref{prop:3}), and we find the~neighbors of
fractions of the~form $\tfrac{1}{j+1}$, $\tfrac{j-1}{2j-1}$, $\tfrac{j}{2j-1}$, $\tfrac{j}{j+1}$
(Corollary~\ref{prop:4}), and of the~form $\tfrac{2}{j+2}$, $\tfrac{j-2}{2(j-1)}$, $\tfrac{j}{2(j-1)}$,
$\tfrac{j}{j+2}$ (Corollary~\ref{prop:6}).

In Section~\ref{section:3}, we slightly simplify our calculations, made in Section~\ref{section:2}, by
describing three subsequences of fractions that are successive in $\mathcal{F}\bigl(\mathbb{B}(2m),m\bigr)$.

In the~paper, we consider Farey (sub)sequences $\mathcal{F}_m$ and $\mathcal{F}\bigl(\mathbb{B}(2m),m\bigr)$
such that $m>1$.

\section{Farey subsequences and collective decision making}
\label{section:applications}

The~Farey sequence $\mathcal{F}_m$ has a~wide area of applications in mathematics, computer science and physics~\cite{CZ}. The~Farey subsequence $\mathcal{F}\bigl(\mathbb{B}(2m),m\bigr)$ which is closely related to $\mathcal{F}_m$, as described in Lemma~\ref{prop:5}, may be useful for analysis of procedures of collective decision making. One such procedure consists in pattern recognition with the~help of committee decision rules, see~\cite{AM-C-I} and references therein.

A~finite collection $\boldsymbol{\mathcal{H}}$ of pairwise distinct hyperplanes in the~{\em feature space\/}~$\mathbb{R}^n$ is a~{\em training set\/} if it is partitioned into two nonempty subsets~$\boldsymbol{\mathcal{A}}$ and~$\boldsymbol{\mathcal{B}}$. A~codimension one subspace~$\pmb{H}:=\{\mathbf{x}\in\mathbb{R}^n:\ \langle\pmb{p},\mathbf{x}\rangle:=\sum_{j=1}^n p_j \mathrm{x}_j=0\}$ from the~{\em hyperplane arrangement\/}~$\boldsymbol{\mathcal{H}}$ is defined by its normal vector $\pmb{p}\in\mathbb{R}^n$, and this hyperplane is oriented: a~vector~$\pmb{v}$ lies on the~{\em positive side\/} of~$\pmb{H}$ if~$\langle\pmb{h},\pmb{v}\rangle>0$, where, by convention, $\pmb{h}:=-\pmb{p}$ if $\pmb{H}\in\boldsymbol{\mathcal{A}}$, and $\pmb{h}:=\pmb{p}$ if $\pmb{H}\in\boldsymbol{\mathcal{B}}$.
In a~similar manner, a~{\em region\/}~$\pmb{T}$ of the~hyperplane arrangement~$\boldsymbol{\mathcal{H}}$, that is, a~connected component of the~{\em complement\/}~$\boldsymbol{\mathcal{T}}:=\mathbb{R}^n-\boldsymbol{\mathcal{H}}$,
lies on the~positive side of the~hyperplane~$\pmb{H}$ if~$\langle\pmb{h},\pmb{v}\rangle>0$ for some vector $\pmb{v}\in\pmb{T}$. Let $\boldsymbol{\mathcal{T}}_{\pmb{H}}^{+}$ denote the~set of all regions lying on the~positive side of $\pmb{H}$.

The~oriented hyperplanes from the~arrangement~$\boldsymbol{\mathcal{H}}$ are called the~{\em training patterns}. The~{\em training samples\/}~$\boldsymbol{\mathcal{A}}$ and $\boldsymbol{\mathcal{B}}$ provide
a~partial description of two disjont {\em classes\/}~$\mathbf{A}$ and~$\mathbf{B}$, respectively: {\sl
a~priori}, we have $\mathbf{A}\supseteq\boldsymbol{\mathcal{A}}$ and~$\mathbf{B}\supseteq\boldsymbol{\mathcal{B}}$.

A~subset $\boldsymbol{\mathcal{K}}^{\ast}\subset\boldsymbol{\mathcal{T}}$ is a~{\em committee of regions\/} for the~arrangement $\boldsymbol{\mathcal{H}}$ if~$|\boldsymbol{\mathcal{K}}^{\ast}\cap\boldsymbol{\mathcal{T}}_{\pmb{H}}^{+}|>\tfrac{1}{2}|\boldsymbol{\mathcal{K}}^{\ast}|$, for each hyperplane $\pmb{H}\in\boldsymbol{\mathcal{H}}$.

Consider a~new pattern $\pmb{G}:=\{\mathbf{x}\in\mathbb{R}^n:\ \langle\pmb{g},\mathbf{x}\rangle=0\}\not\in\boldsymbol{\mathcal{H}}$, defined by its normal vector~$\pmb{g}\in\mathbb{R}^n$. If a~system of distinct representatives $\pmb{W}:=\{\pmb{w}\in\pmb{K}:\ \pmb{K}\in\boldsymbol{\mathcal{K}}^{\ast}\}$, of cardinality $|\boldsymbol{\mathcal{K}}^{\ast}|$, for the~committee of regions $\boldsymbol{\mathcal{K}}^{\ast}$ is fixed, then the~corresponding {\em committee decision rule\/} recognizes the~pattern $\pmb{G}$ as an~element of the~class~$\mathbf{A}$ if $|\{\pmb{w}\in\pmb{W}:\ \langle\pmb{g},\pmb{w}\rangle>0\}|<\tfrac{1}{2}|\pmb{W}|$; the~pattern $\pmb{G}$ is recognized as an~element of the~class $\mathbf{B}$ if $|\{\pmb{w}\in\pmb{W}:\ \langle\pmb{g},\pmb{w}\rangle>0\}|>\tfrac{1}{2}|\pmb{W}|$.

For any hyperplane $\pmb{H}$ from the~arrangement $\boldsymbol{\mathcal{H}}$, the~ascending collection of irreducible fractions
\begin{equation*}
\left(\tfrac{|\boldsymbol{\mathcal{R}}\cap
\boldsymbol{\mathcal{T}}_{\pmb{H}}^{+}|}{\gcd(|\boldsymbol{\mathcal{R}}\cap \boldsymbol{\mathcal{T}}_{\pmb{H}}^{+}|,|\boldsymbol{\mathcal{R}}|)}\!\!\Bigm/\!\!\!
\tfrac{|\boldsymbol{\mathcal{R}}|}{\gcd(|\boldsymbol{\mathcal{R}}\cap
\boldsymbol{\mathcal{T}}_{\pmb{H}}^{+}|,|\boldsymbol{\mathcal{R}}|)}:\ \boldsymbol{\mathcal{R}}\subseteq \boldsymbol{\mathcal{T}},\
|\boldsymbol{\mathcal{R}}|>0\right)
\end{equation*}
is the~Farey subsequence $\mathcal{F}\!\bigl(\mathbb{B}(|\boldsymbol{\mathcal{T}}|),\tfrac{|\boldsymbol{\mathcal{T}}|}{2}\bigr)$, a~dual of the~standard Farey sequence $\mathcal{F}_{|\boldsymbol{\mathcal{T}}|/2}$. A~neighborhood of the~critical value $\tfrac{1}{2}$ in $\mathcal{F}\!\bigl(\mathbb{B}(|\boldsymbol{\mathcal{T}}|),\tfrac{|\boldsymbol{\mathcal{T}}|}{2}\bigr)$ has the~simple structure explained in~Remark~\ref{r}(ii)(iii) of Section~\ref{section:3}.
From the~number-theoretic point of view, a~subset of regions $\boldsymbol{\mathcal{K}}^{\ast}\subset\boldsymbol{\mathcal{T}}$ is a~committee for the~arrangement $\boldsymbol{\mathcal{H}}$ if and only if for each hyperplane $\pmb{H}\in\boldsymbol{\mathcal{H}}$ it holds
\begin{equation*}
\tfrac{|\boldsymbol{\mathcal{K}}^{\ast}\cap
\boldsymbol{\mathcal{T}}_{\pmb{H}}^{+}|}{\gcd(|\boldsymbol{\mathcal{K}}^{\ast}\cap \boldsymbol{\mathcal{T}}_{\pmb{H}}^{+}|,|\boldsymbol{\mathcal{K}}^{\ast}|)}\!\!\Bigm/\!\!\!
\tfrac{|\boldsymbol{\mathcal{K}}^{\ast}|}{\gcd(|\boldsymbol{\mathcal{K}}^{\ast}\cap
\boldsymbol{\mathcal{T}}_{\pmb{H}}^{+}|,|\boldsymbol{\mathcal{K}}^{\ast}|)}
\in\mathcal{F}^{\geq\frac{1}{2}}\!\bigl(\mathbb{B}(|\boldsymbol{\mathcal{T}}|),
\tfrac{|\boldsymbol{\mathcal{T}}|}{2}\bigr)-\{\tfrac{1}{2}\}\ ;
\end{equation*}
thus, the~study of the~structure of the~family of all committees of regions for the~hyperplane arrangement $\boldsymbol{\mathcal{H}}$ might involve the~sequence~$\mathcal{F}^{\geq\frac{1}{2}}\!\bigl(\mathbb{B}(|\boldsymbol{\mathcal{T}}|),
\tfrac{|\boldsymbol{\mathcal{T}}|}{2}\bigr)$.
Such an~attempt is made in~\cite{AM-C-II} within the~bounds of {\em oriented matroid\/} theory.

\section{Neighbors in $\mathcal{F}_m$}
\label{section:1}

We begin this section by recalling several observations made in elementary number theory, cf.~\cite[Chapter~27]{B},
\cite[Chapter~III]{HW}:

$\bullet$ Consider a~fraction $\tfrac{h}{k}\in\mathcal{F}_m-\{\tfrac{0}{1}\}$. To find the~fraction that
precedes~$\tfrac{h}{k}$, consider any integer $x_0$ such that
\begin{equation*}
kx_0\equiv-1\pmod{h}\ ,
\end{equation*}
and $y_0:=\tfrac{kx_0+1}{h}$. For any integer $t$, the~pair $(x_0+ht,y_0+kt)$ is a~solution to the~Diophantine equation $h\mathrm{y}-k\mathrm{x}=1$. The~integer
$t^{\ast}:=\left\lfloor\tfrac{hm-kx_0-1}{hk}\right\rfloor$ is the~maximum solution to the~inequality system
\begin{equation}
\label{eq:4} 0\leq x_0+h\mathrm{t}\leq m \ ,\ \ \ 1\leq y_0+k\mathrm{t}\leq m\ ;
\end{equation}
we have $\tfrac{x_0+ht^{\ast}}{y_0+kt^{\ast}}\in\mathcal{F}_m$. Notice that for any integer solutions $t'$
and $t''$ to system~(\ref{eq:4}), such that $t'\leq t''$, it holds
$\tfrac{0}{1}\leq\tfrac{x_0+ht'}{y_0+kt'}\leq\tfrac{x_0+ht''}{y_0+kt''}<\tfrac{h}{k}$. Since there is no
fraction $\tfrac{i}{j}\in\mathcal{F}_m$ such that
$\tfrac{x_0+ht^{\ast}}{y_0+kt^{\ast}}<\tfrac{i}{j}<\tfrac{h}{k}$, the~fraction
\begin{equation}
\label{eq:5} \left(x_0+h\left\lfloor\tfrac{hm-kx_0-1}{hk}\right\rfloor\right)
\!\Bigm/\!\!\left(\tfrac{kx_0+1}{h}+k\left\lfloor\tfrac{hm-kx_0-1}{hk}\right\rfloor\right)
\end{equation}
precedes $\tfrac{h}{k}$ in~$\mathcal{F}_m$. The~multiplier $\left\lfloor\tfrac{hm-kx_0-1}{hk}\right\rfloor$
in formula~(\ref{eq:5}) turns into zero whenever $\left\lceil\tfrac{hm}{k}\right\rceil-h\leq x_0\leq
\left\lceil\tfrac{hm}{k}\right\rceil-1$.

Now, if $y_0$ is an~integer such that
\begin{equation*}
hy_0\equiv 1\pmod{k}\ ,
\end{equation*}
then the~fraction
\begin{equation}
\label{eq:6} \left(\tfrac{hy_0-1}{k}+h\left\lfloor\tfrac{m-y_0}{k}\right\rfloor\right)
\!\Bigm/\!\!\left(y_0+k\left\lfloor\tfrac{m-y_0}{k}\right\rfloor\right)\ ,
\end{equation}
which coincides with the~fraction described by formula~(\ref{eq:5}), precedes $\tfrac{h}{k}$ in~$\mathcal{F}_m$.

$\bullet$ Consider a~fraction $\tfrac{h}{k}\in\mathcal{F}_m-\{\tfrac{1}{1}\}$.

If $x_0$ and $y_0$ are integers such that
\begin{align*}
kx_0&\equiv 1\pmod{h}\ ,\\ hy_0&\equiv -1\pmod{k}\ ,
\end{align*}
then the~fraction
\begin{align}
\label{eq:7} \left(x_0+h\left\lfloor\tfrac{hm-kx_0+1}{hk}\right\rfloor\right)
\!&\Bigm/\!\!\left(\tfrac{kx_0-1}{h}+k\left\lfloor\tfrac{hm-kx_0+1}{hk}\right\rfloor\right)\\ \label{eq:8} =
\left(\tfrac{hy_0+1}{k}+h\left\lfloor\tfrac{m-y_0}{k}\right\rfloor\right)
\!&\Bigm/\!\!\left(y_0+k\left\lfloor\tfrac{m-y_0}{k}\right\rfloor\right)
\end{align}
succeeds $\tfrac{h}{k}$ in $\mathcal{F}_m$.

Recall that $\tfrac{m-1}{m}$ precedes $\tfrac{1}{1}$, and $\tfrac{1}{m}$ succeeds $\tfrac{0}{1}$ in
$\mathcal{F}_m$. The~following proposition lists simplified versions of expressions~(\ref{eq:5}), (\ref{eq:6}), (\ref{eq:7}) and (\ref{eq:8}) for the~neighbors of fractions in $\mathcal{F}_m$:

\begin{lemma}
\label{prop:1}
\begin{itemize}
\item[\rm(i)] Consider a~fraction $\tfrac{h}{k}\in\mathcal{F}_m-\{\tfrac{0}{1}\}$.
Let $a$ and $b$ be the~integers such that
\begin{align}
\label{eq:3} ka&\equiv -1\pmod{h}\ , & \left\lceil\tfrac{hm}{k}\right\rceil-h\leq\, &a\leq
\left\lceil\tfrac{hm}{k}\right\rceil-1\ ,\\ \nonumber hb&\equiv 1\pmod{k}\ , & m-k+1\leq\, &b\leq m\ .
\end{align}
The~fraction
\begin{equation*}
a\!\!\Bigm/\!\!\tfrac{ka+1}{h}\ \ \ =\ \ \ \tfrac{hb-1}{k}\!\!\Bigm/\!\!b
\end{equation*}
precedes  $\tfrac{h}{k}$ in $\mathcal{F}_m$.

\item[\rm(ii)] Consider a~fraction $\tfrac{h}{k}\in\mathcal{F}_m-\{\tfrac{1}{1}\}$.
Let $a$ and $b$ be the~integers such that
\begin{align}
\label{eq:9} ka&\equiv 1\pmod{h}\ , & \left\lceil\tfrac{hm+2}{k}\right\rceil-h\leq\, &a\leq
\left\lceil\tfrac{hm+2}{k}\right\rceil-1\ ,\\ \nonumber hb&\equiv -1\pmod{k}\ , & m-k+1\leq\, &b\leq m\ .
\end{align}

The~fraction
\begin{equation*}
a\!\!\Bigm/\!\!\tfrac{ka-1}{h}\ \ \ =\ \ \ \tfrac{hb+1}{k}\!\!\Bigm/\!\!b
\end{equation*}
succeeds $\tfrac{h}{k}$ in $\mathcal{F}_m$.
\end{itemize}
\end{lemma}

Inspired by~\cite[Theorem~253]{B}, our approach to the~search for the~neighbors of a~general fraction~$\tfrac{h}{k}$ in~$\mathcal{F}_m$ is more flexible: while the~left boundary of the~search interval in~\cite[Theorem~253]{B} is determined by the~denominator~$k$, the~left boundaries of our additional search intervals included in constraints~(\ref{eq:3}) and~(\ref{eq:9}) depend on the~numerator~$h$; as
a~consequence, in many cases these additional search intervals are much more shorter. Moreover,
if we consider in Lemma~\ref{prop:1} fractions of the~form $\tfrac{1}{j}$ then the~search intervals for
integers~$a$ in constraints~(\ref{eq:3}) and~(\ref{eq:9}) turn into singletons:

\begin{corollary}
\label{prop:2}
\begin{itemize}
\item[\rm(i)]
Consider a~fraction $\tfrac{1}{j}\in\mathcal{F}_m$. The~fraction
\begin{equation}
\label{eq:18} \frac{\left\lceil\frac{m}{j}\right\rceil-1}{j\cdot\left(\left\lceil\frac{m}{j}
\right\rceil-1\right)+1}
\end{equation}
precedes  $\tfrac{1}{j}$ in $\mathcal{F}_m$. If $j>1$, then the~fraction
\begin{equation*}
\frac{\left\lceil\frac{m+2}{j}\right\rceil-1}{j\cdot\left(\left\lceil\frac{m+2}{j} \right\rceil-1\right)-1}
\end{equation*}
succeeds $\tfrac{1}{j}$ in $\mathcal{F}_m$.

\item[\rm(ii)]
Consider a~fraction $\tfrac{j-1}{j}\in\mathcal{F}_m$. If $j>1$, then the~fraction
\begin{equation*}
\frac{(j-1)\left(\left\lceil\frac{m+2}{j}\right\rceil-1\right)-1}{j\cdot\left(\left\lceil\frac{m+2}{j}
\right\rceil-1\right)-1}
\end{equation*}
precedes $\tfrac{j-1}{j}$ in $\mathcal{F}_m$. The~fraction
\begin{equation*}
\frac{(j-1)\left(\left\lceil\frac{m}{j}\right\rceil-1\right)+1}{j\cdot\left(\left\lceil\frac{m}{j}
\right\rceil-1\right)+1}
\end{equation*}
succeeds $\tfrac{j-1}{j}$ in $\mathcal{F}_m$.
\end{itemize}
\end{corollary}

To prove assertion~(ii) of Corollary~\ref{prop:2}, notice that $\tfrac{1}{j}$ is the~image of
$\tfrac{j-1}{j}$ under bijection~(\ref{eq:10}), and find the~images of the~neighbors of $\tfrac{1}{j}$ under
map~(\ref{eq:10}).

The~following statement summarizes the~results of similar calculations that can be performed, with the~help
of Lemma~\ref{prop:1}, for fractions of the~form $\tfrac{2}{j}$, $\tfrac{j-2}{j}\in\mathcal{F}_m$:

\begin{corollary}
\label{prop:7}
\begin{itemize}
\item[\rm(i)] If $\tfrac{2}{j}\in\mathcal{F}_m$, for some $j$, then the~fraction
\begin{equation*}
\begin{cases}\left(\left\lceil\tfrac{2m}{j}\right\rceil-1\right)\!\!\biggm/
\!\!\tfrac{j\cdot\left(\left\lceil\frac{2m}{j}\right\rceil -1\right)+1}{2}\ ,& \text{if
$\left\lceil\tfrac{2m}{j}\right\rceil\equiv 0\pmod{2}$}\ ;\\
\left(\left\lceil\tfrac{2m}{j}\right\rceil-2\right)\!\!\biggm/
\!\!\tfrac{j\cdot\left(\left\lceil\frac{2m}{j}\right\rceil-2\right)+1}{2}\ ,& \text{if
$\left\lceil\tfrac{2m}{j}\right\rceil\equiv 1\pmod{2}$}\ ;
\end{cases}
\end{equation*}
precedes $\tfrac{2}{j}$ in $\mathcal{F}_m$; the~fraction
\begin{equation*}
\begin{cases}\left(\left\lceil\tfrac{2(m+1)}{j}\right\rceil-1\right)\!\!\biggm/
\!\!\tfrac{j\cdot\left(\left\lceil\frac{2(m+1)}{j}\right\rceil-1\right)-1}{2}\ ,& \text{if
$\left\lceil\tfrac{2(m+1)}{j}\right\rceil\equiv 0\pmod{2}$}\ ;\\
\left(\left\lceil\tfrac{2(m+1)}{j}\right\rceil-2\right)\!\!\biggm/
\!\!\tfrac{j\cdot\left(\left\lceil\frac{2(m+1)}{j}\right\rceil-2\right)-1}{2}\ ,& \text{if
$\left\lceil\tfrac{2(m+1)}{j}\right\rceil\equiv 1\pmod{2}$}\ ;
\end{cases}
\end{equation*}
succeeds $\tfrac{2}{j}$ in $\mathcal{F}_m$.

\item[\rm(ii)] If $\tfrac{j-2}{j}\in\mathcal{F}_m$, for some $j$, then the~fraction
\begin{equation*}
\begin{cases}\tfrac{(j-2)\left(\left\lceil\frac{2(m+1)}{j}\right\rceil-1\right)-1}{2}\!\!\biggm/
\!\!\tfrac{j\cdot\left(\left\lceil\frac{2(m+1)}{j}\right\rceil-1\right)-1}{2}\ ,& \text{if
$\left\lceil\tfrac{2(m+1)}{j}\right\rceil\equiv 0\pmod{2}$}\ ;\\
\tfrac{(j-2)\left(\left\lceil\frac{2(m+1)}{j}\right\rceil-2\right)-1}{2}\!\!\biggm/
\!\!\tfrac{j\cdot\left(\left\lceil\frac{2(m+1)}{j}\right\rceil-2\right)-1}{2}\ ,& \text{if
$\left\lceil\tfrac{2(m+1)}{j}\right\rceil\equiv 1\pmod{2}$}\ ;
\end{cases}
\end{equation*}
precedes $\tfrac{j-2}{j}$ in $\mathcal{F}_m$; the~fraction
\begin{equation*}
\begin{cases}\tfrac{(j-2)\left(\left\lceil\frac{2m}{j}\right\rceil -1\right)+1}{2}\!\!\biggm/
\!\!\tfrac{j\cdot\left(\left\lceil\frac{2m}{j}\right\rceil -1\right)+1}{2}\ ,& \text{if
$\left\lceil\tfrac{2m}{j}\right\rceil\equiv 0\pmod{2}$}\ ;\\
\tfrac{(j-2)\left(\left\lceil\frac{2m}{j}\right\rceil-2\right)+1}{2}\!\!\biggm/
\!\!\tfrac{j\cdot\left(\left\lceil\frac{2m}{j}\right\rceil-2\right)+1}{2}\ ,& \text{if
$\left\lceil\tfrac{2m}{j}\right\rceil\equiv 1\pmod{2}$}\ ;
\end{cases}
\end{equation*}
succeeds $\tfrac{j-2}{j}$ in $\mathcal{F}_m$.
\end{itemize}
\end{corollary}

\section{Neighbors in $\mathcal{F}\bigl(\mathbb{B}(2m),m\bigr)$}
\label{section:2}

We now extend the~results, obtained in the~previous section, to the~Farey subsequence
$\mathcal{F}\bigl(\mathbb{B}(2m),m\bigr)$. We begin by searching for the~neighbors of an~arbitrary fraction
in $\mathcal{F}\bigl(\mathbb{B}(2m),m\bigr)$:

\begin{proposition}
\label{prop:3}
\begin{itemize}
\item[\rm(i)]
Consider a~fraction $\tfrac{h}{k}\in\mathcal{F}^{\leq\frac{1}{2}}\bigl(\mathbb{B}(2m),m\bigr)$.
\begin{itemize}
\item[\rm(a)] If $\tfrac{h}{k}\neq\tfrac{0}{1}$, let $a$ and $b$ be the~integers such that
\begin{align*}
(k-h)a&\equiv -1\pmod{h}\ , & \left\lceil\tfrac{hm}{k-h}\right\rceil-h\leq\, &a\leq
\left\lceil\tfrac{hm}{k-h}\right\rceil-1\ ,\\ hb&\equiv 1\pmod{(k-h)}\ , & m-k+h+1\leq\, &b\leq m\ .
\end{align*}
The~fraction
\begin{equation*}
a\!\!\Bigm/\!\!\tfrac{ka+1}{h}\ \ \ =\ \ \ \tfrac{hb-1}{k-h}\!\!\Bigm/\!\!\tfrac{kb-1}{k-h}
\end{equation*}
precedes  $\tfrac{h}{k}$ in $\mathcal{F}\bigl(\mathbb{B}(2m),m\bigr)$.

\item[\rm(b)]
If $\tfrac{h}{k}\neq\tfrac{1}{2}$, let $a$ and $b$ be the~integers such that
\begin{align*}
(k-h)a&\equiv 1\pmod{h}\ , & \left\lceil\tfrac{hm+2}{k-h}\right\rceil-h\leq\, &a\leq
\left\lceil\tfrac{hm+2}{k-h}\right\rceil-1\ ,\\ hb&\equiv -1\pmod{(k-h)}\ , & m-k+h+1\leq\, &b\leq m\ .
\end{align*}

The~fraction
\begin{equation*}
a\!\!\Bigm/\!\!\tfrac{ka-1}{h}\ \ \ =\ \ \ \tfrac{hb+1}{k-h}\!\!\Bigm/\!\!\tfrac{kb+1}{k-h}
\end{equation*}
succeeds $\tfrac{h}{k}$ in $\mathcal{F}\bigl(\mathbb{B}(2m),m\bigr)$.
\end{itemize}

\item[\rm(ii)]
Consider a~fraction $\tfrac{h}{k}\in\mathcal{F}^{\geq\frac{1}{2}}\bigl(\mathbb{B}(2m),m\bigr)$.

\begin{itemize}
\item[\rm(a)] If $\tfrac{h}{k}\neq\tfrac{1}{2}$, let $a$ and $b$ be the~integers such that
\begin{align*}
ka&\equiv -1\pmod{h}\ , & m-h+1\leq\, &a\leq m\ ,\\ hb&\equiv 1\pmod{(k-h)}\ , &
\left\lceil\tfrac{(k-h)m+2}{h}\right\rceil-k+h\leq\, &b\leq \left\lceil\tfrac{(k-h)m+2}{h}\right\rceil-1\ .
\end{align*}
The~fraction
\begin{align*}
a\!\!\Bigm/\!\!\tfrac{ka+1}{h}\ \ \ =\ \ \ \tfrac{hb-1}{k-h}\!\!\Bigm/\!\!\tfrac{kb-1}{k-h}
\end{align*}
precedes  $\tfrac{h}{k}$ in $\mathcal{F}\bigl(\mathbb{B}(2m),m\bigr)$.

\item[\rm(b)]
If $\tfrac{h}{k}\neq\tfrac{1}{1}$, let $a$ and $b$ be the~integers such that
\begin{align*}
ka&\equiv 1\pmod{h}\ , & m-h+1\leq\, &a\leq m\ ,\\ hb&\equiv -1\pmod{(k-h)}\ , &
\left\lceil\tfrac{(k-h)m}{h}\right\rceil-k+h\leq\, &b\leq \left\lceil\tfrac{(k-h)m}{h}\right\rceil-1\ .
\end{align*}

The~fraction
\begin{equation*}
a\!\!\Bigm/\!\!\tfrac{ka-1}{h}\ \ \ =\ \ \ \tfrac{hb+1}{k-h}\!\!\Bigm/\!\!\tfrac{kb+1}{k-h}
\end{equation*}
succeeds $\tfrac{h}{k}$ in $\mathcal{F}\bigl(\mathbb{B}(2m),m\bigr)$.
\end{itemize}
\end{itemize}
\end{proposition}

\begin{proof} To prove assertion~(i)(a), notice that $\tfrac{h}{k-h}\in\mathcal{F}_m$ is the~image
of $\tfrac{h}{k}$ under bijection~(\ref{eq:11}), use~Lemma~\ref{prop:1}(i) to find the~predecessor of
$\tfrac{h}{k-h}$ in $\mathcal{F}_m$, and send it to~$\mathcal{F}\bigl(\mathbb{B}(2m),m\bigr)$ by means of
bijection~(\ref{eq:12}).

Assertion~(i)(b)\hfill is\hfill proved\hfill in\hfill a\hfill similar\hfill way,\hfill by\hfill the~application\\ of~Lemma~\ref{prop:1}(ii), with the~help of
bijections~(\ref{eq:11}) and~(\ref{eq:12}).

One can prove Proposition~\ref{prop:3}(ii) by the~application of~Lemma~\ref{prop:1}, with the~help of
bijections~(\ref{eq:13}) and~(\ref{eq:14}).
\end{proof}

\begin{corollary}
\label{prop:4}
\begin{itemize}
\item[\rm(i)]
If $\tfrac{1}{j+1}\in\mathcal{F}\bigl(\mathbb{B}(2m),m\bigr)$ and $\tfrac{1}{j+1}<\tfrac{1}{2}$, for some
$j$, then the~fraction
\begin{equation*}
\frac{\left\lceil\frac{m}{j}\right\rceil-1}{(j+1)\left(\left\lceil\frac{m}{j} \right\rceil-1\right)+1}
\end{equation*}
precedes  $\tfrac{1}{j+1}$ in $\mathcal{F}\bigl(\mathbb{B}(2m),m\bigr)$; the~fraction
\begin{equation*}
\frac{\left\lceil\frac{m+2}{j}\right\rceil-1}{(j+1)\left(\left\lceil\frac{m+2}{j} \right\rceil-1\right)-1}
\end{equation*}
succeeds $\tfrac{1}{j+1}$ in $\mathcal{F}\bigl(\mathbb{B}(2m),m\bigr)$.

\item[\rm(ii)]
If $\frac{j-1}{2j-1}\in\mathcal{F}\bigl(\mathbb{B}(2m),m\bigr)$ and $\frac{j-1}{2j-1}<\frac{1}{2}$, for some
$j$, then the~fraction
\begin{equation*}
\frac{(j-1)\left(\left\lceil\frac{m+2}{j}\right\rceil-1\right)-1}{(2j-1)\left(\left\lceil\frac{m+2}{j}
\right\rceil-1\right)-2}
\end{equation*}
precedes $\frac{j-1}{2j-1}$ in $\mathcal{F}\bigl(\mathbb{B}(2m),m\bigr)$; the~fraction
\begin{equation*}
\frac{(j-1)\left(\left\lceil\frac{m}{j}\right\rceil-1\right)+1}{(2j-1)\left(\left\lceil\frac{m}{j}
\right\rceil-1\right)+2}
\end{equation*}
succeeds $\frac{j-1}{2j-1}$ in $\mathcal{F}\bigl(\mathbb{B}(2m),m\bigr)$.

\item[\rm(iii)]
If $\frac{j}{2j-1}\in\mathcal{F}\bigl(\mathbb{B}(2m),m\bigr)$ and $\frac{1}{2}<\frac{j}{2j-1}<\tfrac{1}{1}$,
for some $j$, then the~fraction
\begin{equation*}
\frac{j\cdot\left(\left\lceil\frac{m}{j}\right\rceil-1\right)+1}{(2j-1)\left(\left\lceil\frac{m}{j}
\right\rceil-1\right)+2}
\end{equation*}
precedes $\frac{j}{2j-1}$ in $\mathcal{F}\bigl(\mathbb{B}(2m),m\bigr)$; the~fraction
\begin{equation*}
\frac{j\cdot\left(\left\lceil\frac{m+2}{j}\right\rceil-1\right)-1}{(2j-1)\left(\left\lceil\frac{m+2}{j}
\right\rceil-1\right)-2}
\end{equation*}
succeeds $\frac{j}{2j-1}$ in $\mathcal{F}\bigl(\mathbb{B}(2m),m\bigr)$.

\item[\rm(iv)]
If $\frac{j}{j+1}\in\mathcal{F}\bigl(\mathbb{B}(2m),m\bigr)$ and $\frac{1}{2}<\frac{j}{j+1}$, for some $j$,
then the~fraction
\begin{equation*}
\frac{j\cdot\left(\left\lceil\frac{m+2}{j}\right\rceil-1\right)-1}{(j+1)\left(\left\lceil\frac{m+2}{j}
\right\rceil-1\right)-1}
\end{equation*}
precedes $\frac{j}{j+1}$ in $\mathcal{F}\bigl(\mathbb{B}(2m),m\bigr)$; the~fraction
\begin{equation*}
\frac{j\cdot\left(\left\lceil\frac{m}{j}\right\rceil-1\right)+1}{(j+1)\left(\left\lceil\frac{m}{j}
\right\rceil-1\right)+1}
\end{equation*}
succeeds $\frac{j}{j+1}$ in $\mathcal{F}\bigl(\mathbb{B}(2m),m\bigr)$.
\end{itemize}
\end{corollary}

\begin{proof} The~fraction $\tfrac{1}{j}\in\mathcal{F}_m$ is the~image of the~fractions $\tfrac{1}{j+1}$,
$\tfrac{j-1}{2j-1}$, $\tfrac{j}{2j-1}$ and $\tfrac{j}{j+1}$ under bijections~(\ref{eq:11}), (\ref{eq:15}),
(\ref{eq:16}) and (\ref{eq:13}), respectively.

To prove assertion~(i), find the~predecessor of $\tfrac{1}{j}$ in $\mathcal{F}_m$ by the~application of
Corollary~\ref{prop:2}(i), and reflect it to $\mathcal{F}\bigl(\mathbb{B}(2m),m\bigr)$ by means of
bijection~(\ref{eq:12}).

One can prove the~remaining assertions in a~similar way, with the~help of Corollary~\ref{prop:2} and of the~monotone bijections mentioned in Lemma~\ref{prop:5}.
\end{proof}

Recall also that the~fraction $\tfrac{m}{m+1}$ precedes $\tfrac{1}{1}$, and $\tfrac{1}{m+1}$ succeeds
$\tfrac{0}{1}$ in $\mathcal{F}\bigl(\mathbb{B}(2m),m\bigr)$; the~three fractions
$\tfrac{m-1}{2m-1}<\tfrac{1}{2}<\tfrac{m}{2m-1}$ are successive in $\mathcal{F}\bigl(\mathbb{B}(2m),m\bigr)$,
see~\cite{AM-C-II}.

\begin{corollary}
\label{prop:6}
\begin{itemize}
\item[\rm(i)]
If $\tfrac{2}{j+2}\in\mathcal{F}\bigl(\mathbb{B}(2m),m\bigr)$ and $\tfrac{2}{j+2}<\tfrac{1}{2}$, for some
$j$, then the~fraction
\begin{equation*}
\begin{cases}\left(\left\lceil\tfrac{2m}{j}\right\rceil-1\right)\!\!\biggm/
\!\!\tfrac{(j+2)\left(\left\lceil\frac{2m}{j}\right\rceil -1\right)+1}{2}\ ,& \text{if
$\left\lceil\tfrac{2m}{j}\right\rceil\equiv 0\pmod{2}$}\ ;\\
\left(\left\lceil\tfrac{2m}{j}\right\rceil-2\right)\!\!\biggm/
\!\!\tfrac{(j+2)\left(\left\lceil\frac{2m}{j}\right\rceil-2\right)+1}{2}\ ,& \text{if
$\left\lceil\tfrac{2m}{j}\right\rceil\equiv 1\pmod{2}$}\ ;
\end{cases}
\end{equation*}
precedes $\tfrac{2}{j+2}$ in $\mathcal{F}\bigl(\mathbb{B}(2m),m\bigr)$; the~fraction
\begin{equation*}
\begin{cases}\left(\left\lceil\tfrac{2(m+1)}{j}\right\rceil-1\right)\!\!\biggm/
\!\!\tfrac{(j+2)\left(\left\lceil\frac{2(m+1)}{j}\right\rceil-1\right)-1}{2}\ ,& \text{if
$\left\lceil\tfrac{2(m+1)}{j}\right\rceil\equiv 0\pmod{2}$}\ ;\\
\left(\left\lceil\tfrac{2(m+1)}{j}\right\rceil-2\right)\!\!\biggm/
\!\!\tfrac{(j+2)\left(\left\lceil\frac{2(m+1)}{j}\right\rceil-2\right)-1}{2}\ ,& \text{if
$\left\lceil\tfrac{2(m+1)}{j}\right\rceil\equiv 1\pmod{2}$}\ ;
\end{cases}
\end{equation*}
succeeds $\tfrac{2}{j+2}$ in $\mathcal{F}\bigl(\mathbb{B}(2m),m\bigr)$.

\item[\rm(ii)]
If $\tfrac{j-2}{2(j-1)}\in\mathcal{F}\bigl(\mathbb{B}(2m),m\bigr)$ and $\tfrac{j-2}{2(j-1)}<\tfrac{1}{2}$,
for some $j$, then the~fraction
\begin{equation*}
\begin{cases}\tfrac{(j-2)\left(\left\lceil\frac{2(m+1)}{j}\right\rceil-1\right)-1}{2}\!\!\biggm/\!\!
\left((j-1)\left(\left\lceil\frac{2(m+1)}{j}\right\rceil-1\right)-1\right)\ ,& \text{if
$\left\lceil\tfrac{2(m+1)}{j}\right\rceil\equiv 0\pmod{2}$}\ ;\\
\tfrac{(j-2)\left(\left\lceil\frac{2(m+1)}{j}\right\rceil-2\right)-1}{2}\!\!\biggm/ \!\!
\left((j-1)\left(\left\lceil\tfrac{2(m+1)}{j}\right\rceil-2\right)-1\right)\ ,& \text{if
$\left\lceil\tfrac{2(m+1)}{j}\right\rceil\equiv 1\pmod{2}$}\ ;
\end{cases}
\end{equation*}
precedes $\tfrac{j-2}{2(j-1)}$ in $\mathcal{F}\bigl(\mathbb{B}(2m),m\bigr)$; the~fraction
\begin{equation*}
\begin{cases}
\tfrac{(j-2)\left(\left\lceil\frac{2m}{j}\right\rceil-1\right)+1}{2}\!\!\biggm/ \!\!
\left((j-1)\left(\left\lceil\tfrac{2m}{j}\right\rceil-1\right)+1\right)\ ,& \text{if
$\left\lceil\tfrac{2m}{j}\right\rceil\equiv 0\pmod{2}$}\ ;\\
\tfrac{(j-2)\left(\left\lceil\frac{2m}{j}\right\rceil-2\right)+1}{2}\!\!\biggm/\!\!
\left((j-1)\left(\left\lceil\tfrac{2m}{j}\right\rceil-2\right)+1\right)\ ,& \text{if
$\left\lceil\tfrac{2m}{j}\right\rceil\equiv 1\pmod{2}$}\ ;
\end{cases}
\end{equation*}
succeeds $\tfrac{j-2}{2(j-1)}$ in $\mathcal{F}\bigl(\mathbb{B}(2m),m\bigr)$.

\item[\rm(iii)]
If $\tfrac{j}{2(j-1)}\in\mathcal{F}\bigl(\mathbb{B}(2m),m\bigr)$ and $\tfrac{1}{2}<\tfrac{j}{2(j-1)}$, for
some $j$, then the~fraction
\begin{equation*}
\begin{cases}
\tfrac{j\cdot\left(\left\lceil\frac{2m}{j}\right\rceil -1\right)+1}{2} \!\!\biggm/ \!\!
\left((j-1)\left(\left\lceil\tfrac{2m}{j}\right\rceil-1\right)+1\right)\ ,& \text{if
$\left\lceil\tfrac{2m}{j}\right\rceil\equiv 0\pmod{2}$}\ ;\\
\tfrac{j\cdot\left(\left\lceil\frac{2m}{j}\right\rceil-2\right)+1}{2}\!\!\biggm/\!\!
\left((j-1)\left(\left\lceil\tfrac{2m}{j}\right\rceil-2\right)+1\right)\ ,& \text{if
$\left\lceil\tfrac{2m}{j}\right\rceil\equiv 1\pmod{2}$}\ ;
\end{cases}
\end{equation*}
precedes $\tfrac{j}{2(j-1)}$ in $\mathcal{F}\bigl(\mathbb{B}(2m),m\bigr)$; the~fraction
\begin{equation*}
\begin{cases}
\tfrac{j\cdot\left(\left\lceil\frac{2(m+1)}{j}\right\rceil-1\right)-1}{2} \!\!\biggm/ \!\!
\left((j-1)\left(\left\lceil\tfrac{2(m+1)}{j}\right\rceil-1\right)-1\right)\ ,& \text{if
$\left\lceil\tfrac{2(m+1)}{j}\right\rceil\equiv 0\pmod{2}$}\ ;\\
\tfrac{j\cdot\left(\left\lceil\frac{2(m+1)}{j}\right\rceil-2\right)-1}{2} \!\!\biggm/\!\!
\left((j-1)\left(\left\lceil\tfrac{2(m+1)}{j}\right\rceil-2\right)-1\right)\ ,& \text{if
$\left\lceil\tfrac{2(m+1)}{j}\right\rceil\equiv 1\pmod{2}$}\ ;
\end{cases}
\end{equation*}
succeeds $\tfrac{j}{2(j-1)}$ in $\mathcal{F}\bigl(\mathbb{B}(2m),m\bigr)$.

\item[\rm(iv)]
If $\tfrac{j}{j+2}\in\mathcal{F}\bigl(\mathbb{B}(2m),m\bigr)$ and $\tfrac{1}{2}<\tfrac{j}{j+2}$, for some
$j$, then the~fraction
\begin{equation*}
\begin{cases}\tfrac{j\cdot\left(\left\lceil\frac{2(m+1)}{j}\right\rceil-1\right)-1}{2}
\!\!\biggm/\!\!\tfrac{(j+2)\left(\left\lceil\frac{2(m+1)}{j}\right\rceil-1\right)-1}{2} \ ,& \text{if
$\left\lceil\tfrac{2(m+1)}{j}\right\rceil\equiv 0\pmod{2}$}\ ;\\
\tfrac{j\cdot\left(\left\lceil\frac{2(m+1)}{j}\right\rceil-2\right)-1}{2}\!\!\biggm/
\!\!\tfrac{(j+2)\left(\left\lceil\frac{2(m+1)}{j}\right\rceil-2\right)-1}{2}\ ,& \text{if
$\left\lceil\tfrac{2(m+1)}{j}\right\rceil\equiv 1\pmod{2}$}\ ;
\end{cases}
\end{equation*}
precedes $\tfrac{j}{j+2}$ in $\mathcal{F}\bigl(\mathbb{B}(2m),m\bigr)$; the~fraction
\begin{equation*}
\begin{cases}\tfrac{j\cdot\left(\left\lceil\frac{2m}{j}\right\rceil-1\right)+1}{2}\!\!\biggm/
\!\!\tfrac{(j+2)\left(\left\lceil\frac{2m}{j}\right\rceil -1\right)+1}{2}\ ,& \text{if
$\left\lceil\tfrac{2m}{j}\right\rceil\equiv 0\pmod{2}$}\ ;\\
\tfrac{j\cdot\left(\left\lceil\frac{2m}{j}\right\rceil-2\right)+1}{2}\!\!\biggm/
\!\!\tfrac{(j+2)\left(\left\lceil\frac{2m}{j}\right\rceil-2\right)+1}{2}\ ,& \text{if
$\left\lceil\tfrac{2m}{j}\right\rceil\equiv 1\pmod{2}$}\ ;
\end{cases}
\end{equation*}
succeeds $\tfrac{j}{j+2}$ in $\mathcal{F}\bigl(\mathbb{B}(2m),m\bigr)$.
\end{itemize}
\end{corollary}

\begin{proof} Assertions~(i) and~(ii) are reformulations of Corollary~\ref{prop:7}(i), made with the~help of
monotone bijections~(\ref{eq:11}), (\ref{eq:12}), (\ref{eq:15}) and (\ref{eq:21}). Corollary~\ref{prop:7}(ii)
leads to assertions~(iii) and~(iv) by means of bijections~(\ref{eq:13}), (\ref{eq:14}), (\ref{eq:16}) and~(\ref{eq:20}).
\end{proof}

\section{Three subsequences of adjacent fractions within~$\mathcal{F}\bigl(\mathbb{B}(2m),m\bigr)$ }
\label{section:3}

Formula~(\ref{eq:18}) implies that the~fractions
\begin{equation}
\label{eq:19}
\tfrac{0}{1}<\tfrac{1}{m}<\tfrac{1}{m-1}<\tfrac{1}{m-2}<\cdots<\tfrac{1}{\lceil m/2\rceil}
\end{equation}
are consecutive in $\mathcal{F}_m$, in the~same way as the~fractions
\begin{equation*}
\tfrac{\lceil m/2\rceil-1}{\lceil
m/2\rceil}<\cdots<\tfrac{m-3}{m-2}<\tfrac{m-2}{m-1}<\tfrac{m-1}{m}<\tfrac{1}{1}
\end{equation*}
are consecutive in $\mathcal{F}_m$ thanks to bijection~(\ref{eq:10}); therefore we can clarify the~statement
of Corollary~\ref{prop:4} in the~following way:

\begin{remark}
\label{r}
\begin{itemize}
\item[\rm(i)]
The~fractions
\begin{equation*}
\tfrac{0}{1}<\tfrac{1}{m+1}<\tfrac{1}{m}<\tfrac{1}{m-1}<\cdots<\tfrac{1}{\lceil m/2\rceil+1}
\end{equation*}
are consecutive in $\mathcal{F}^{\leq\frac{1}{2}}\bigl(\mathbb{B}(2m),m\bigr)$.

\item[\rm(ii)]
The~fractions
\begin{equation*}
\tfrac{\lceil m/2\rceil-1}{2\lceil
m/2\rceil-1}<\cdots<\tfrac{m-3}{2m-5}<\tfrac{m-2}{2m-3}<\tfrac{m-1}{2m-1}<\tfrac{1}{2}
\end{equation*}
are consecutive in $\mathcal{F}^{\leq\frac{1}{2}}\bigl(\mathbb{B}(2m),m\bigr)$.

\item[\rm(iii)]
The~fractions
\begin{equation*}
\tfrac{1}{2}<\tfrac{m}{2m-1}<\tfrac{m-1}{2m-3}<\tfrac{m-2}{2m-5}<\cdots<\tfrac{\lceil m/2\rceil}{2\lceil
m/2\rceil-1}
\end{equation*}
are consecutive in $\mathcal{F}^{\geq\frac{1}{2}}\bigl(\mathbb{B}(2m),m\bigr)$.

\item[\rm(iv)]
The~fractions
\begin{equation*}
\tfrac{\lceil m/2\rceil}{\lceil
m/2\rceil+1}<\cdots<\tfrac{m-2}{m-1}<\tfrac{m-1}{m}<\tfrac{m}{m+1}<\tfrac{1}{1}
\end{equation*}
are consecutive in $\mathcal{F}^{\geq\frac{1}{2}}\bigl(\mathbb{B}(2m),m\bigr)$.
\end{itemize}
Indeed, since the~fractions composing sequence~{\rm(\ref{eq:19})} are consecutive in $\mathcal{F}_m$, we
arrive at conclusions~{\rm(i)}, {\rm(ii)}, {\rm(iii)} and {\rm(iv)} with the~help of monotone
bijections~{\rm(\ref{eq:12})}, {\rm(\ref{eq:21})}, {\rm(\ref{eq:20})} and {\rm(\ref{eq:14})}, respectively.
\end{remark}

\end{document}